\documentclass[submission,copyright,creativecommons]{eptcs}

\newcommand{\onetitle}{DisCoPy: the Hierarchy of Graphical Languages in Python}

\providecommand{\event}{ACT 2023} 

\usepackage{iftex}
\usepackage[english]{babel}

\ifpdf
  \usepackage{underscore}         
  \usepackage[T1]{fontenc}        

\else
  \usepackage{breakurl}           
\fi

\title{\onetitle}
\author{Alexis Toumi \qquad\qquad Richie Yeung \qquad\qquad Boldizs{\'a}r Po{\'o}r  \qquad\qquad Giovanni de Felice
\institute{Quantinuum – Quantum Compositional Intelligence\\
17 Beaumont street, OX1 2NA Oxford, UK}
\email{firstname@discopy.org}}

\newcommand{\titlerunning}{\onetitle}
\newcommand{\authorrunning}{A. Toumi, R. Yeung, B. Po{\'o}r \& G. de Felice}

\hypersetup{
  bookmarksnumbered,
  pdftitle    = {\titlerunning},
  pdfauthor   = {\authorrunning},
  pdfsubject  = {EPTCS},               
}

\usepackage[frozencache]{minted}
\setminted{frame=lines,fontseries=mono,fontsize=\footnotesize}
\setmintedinline{fontsize=auto}

\newcommand{\py}[1]{{\color{brown}\mintinline[breaklines]{python}{#1}}}
\newcommand{\pyref}[2]{\href{https://docs.discopy.org/en/main/_api/discopy.#1.html}{\py{#2}}}

\usepackage{appendix}

\usepackage{xcolor}

\hypersetup{colorlinks=true,urlcolor=orange,linkcolor=orange,citecolor=orange}
\usepackage{tikzit}

\usepackage{caption}
\usepackage{subcaption}

\makeatletter
\newcommand{\ssymbol}[1]{^{\@fnsymbol{#1}}}

\makeatother

\tikzstyle{gate}=[shape=rectangle, text height=1.5ex, text depth=0.25ex, yshift=0.5mm, fill=white, draw=black, minimum height=5mm, yshift=-0.5mm, minimum width=5mm, font={\small}, tikzit category=circuit]
\tikzstyle{big gate}=[shape=rectangle, text height=1.5ex, text depth=0.25ex, yshift=0.5mm, fill=white, draw=black, minimum height=10mm, yshift=-0.5mm, minimum width=5mm, font={\small}, tikzit category=circuit]

\tikzstyle{Z dot}=[inner sep=0mm, minimum size=2mm, shape=circle, draw=black, fill={rgb,255: red,221; green,255; blue,221}, tikzit category=zx]
\tikzstyle{Z phase dot}=[minimum size=1.2em, font={\footnotesize\boldmath}, shape=rectangle, rounded corners=0.5em, inner sep=0.2em, outer sep=-0.2em, scale=0.8, draw=black, fill={rgb,255: red,221; green,255; blue,221}, tikzit shape=circle, tikzit draw=blue, tikzit category=zx]
\tikzstyle{Z box}=[Z phase dot, rounded corners=0, fill={rgb,255: red,221; green,255; blue,221}, tikzit shape=rectangle, tikzit draw=blue, tikzit category=zx]
\tikzstyle{box}=[Z box, fill=white]

\tikzstyle{X dot}=[Z dot, shape=circle, draw=black, fill={rgb,255: red,232; green,165; blue,165}, tikzit category=zx]
\tikzstyle{X phase dot}=[Z phase dot, tikzit shape=circle, tikzit draw=blue, fill={rgb,255: red,232; green,165; blue,165}, font={\footnotesize\boldmath}, tikzit category=zx]

\tikzstyle{red dot}=[Z dot, shape=circle, draw=black, fill={rgb,255: red,255; green,0; blue,0}, tikzit category=zx]
\tikzstyle{red phase dot}=[Z phase dot, text=white, draw=black, fill={rgb,255: red,255; green,0; blue,0}, tikzit draw=blue, tikzit category=zx]

\tikzstyle{hadamard}=[fill=yellow, draw=black, shape=rectangle, inner sep=0.6mm, minimum height=1.5mm, minimum width=1.5mm, tikzit category=zx]
\tikzstyle{dtriangle}=[fill=yellow,draw=black,shape=isosceles triangle,shape border rotate=-90,isosceles triangle stretches=true,inner sep=0.8pt,minimum width=0.25cm,minimum height=2mm]
\tikzstyle{vtriang}=[fill=yellow,draw=black,shape=isosceles triangle,shape border rotate=180,isosceles triangle stretches=true,inner sep=0.8pt,minimum width=0.25cm,minimum height=2mm]
\tikzstyle{bspider}=[fill=black,draw=black,scale=1,shape=isosceles triangle,shape border rotate=-90,isosceles triangle stretches=true,inner sep=1pt,minimum width=0.4cm,minimum height=3mm]
\tikzstyle{dbspider}=[fill=black,draw=black,scale=1,shape=isosceles triangle,shape border rotate=90,isosceles triangle stretches=true,inner sep=1pt,minimum width=0.4cm,minimum height=3mm]

\tikzstyle{starv}=[fill=yellow,draw=black,shape=star,star points=6,star point ratio=1.74,inner sep=0.8pt,minimum height=2.75mm]
\tikzstyle{starh}=[fill=yellow,draw=black,shape=star,star points=6,star point ratio=1.74,inner sep=0.8pt,minimum height=2.75mm, rotate=90]

\tikzstyle{paulibox}=[fill={rgb,255: red,221; green,221; blue,255}, draw=black, shape=rectangle, inner sep=0.6mm, minimum height=5mm, minimum width=5mm, font={\footnotesize}, text height=1.5ex, text depth=0.25ex, tikzit category=zx]
\tikzstyle{vertex}=[inner sep=0mm, minimum size=1mm, shape=circle, draw=black, fill=black, tikzit category=misc]
\tikzstyle{vertex set}=[inner sep=0mm, minimum size=1mm, shape=circle, draw=black, fill=white, font={\footnotesize\boldmath}, tikzit category=misc]
\tikzstyle{small black dot}=[fill=black, draw=black, shape=circle, inner sep=0pt, minimum width=1.2mm, tikzit category=circuit]
\tikzstyle{cnot ctrl}=[fill=black, draw=black, shape=circle, inner sep=0pt, minimum width=1.2mm, tikzit category=circuit]
\tikzstyle{cnot targ}=[fill=white, draw=white, shape=circle, tikzit category=circuit, label={center:$\oplus$}, inner sep=0pt, minimum width=2.1mm, tikzit fill={rgb,255: red,102; green,204; blue,255}, tikzit draw=black]
\tikzstyle{ket}=[fill=white, draw=black, shape=regular polygon, regular polygon sides=3, regular polygon rotate=-30, scale=0.7, inner sep=1pt, tikzit category=circuit, tikzit shape=rectangle, tikzit fill=green]
\tikzstyle{bra}=[fill=white, draw=black, shape=regular polygon, regular polygon sides=3, regular polygon rotate=30, scale=0.7, inner sep=1pt, tikzit category=circuit, tikzit shape=rectangle, tikzit fill=red]
\tikzstyle{scalar}=[shape=rectangle, text height=1.5ex, text depth=0.25ex, yshift=0.5mm, fill=white, draw=black, minimum height=5mm, yshift=-0.5mm, minimum width=5mm, font={\small}]
\tikzstyle{clabel}=[fill=white, draw=none, shape=rectangle, tikzit fill={rgb,255: red,56; green,255; blue,242}, font={\footnotesize}, inner sep=1pt, tikzit category=labels]
\tikzstyle{empty diagram}=[draw={gray!40!white}, dashed, shape=rectangle, minimum width=1cm, minimum height=1cm, tikzit category=misc]
\tikzstyle{amap}=[fill=white, draw=black, shape=NEbox, tikzit category=asymmetric, tikzit fill=yellow, tikzit shape=rectangle]
\tikzstyle{amap conj}=[fill=white, draw=black, shape=NWbox, tikzit category=asymmetric, tikzit fill=green, tikzit shape=rectangle]
\tikzstyle{amap adj}=[fill=white, draw=black, shape=SEbox, tikzit category=asymmetric, tikzit fill=red, tikzit shape=rectangle]
\tikzstyle{amap trans}=[fill=white, draw=black, shape=SWbox, tikzit category=asymmetric, tikzit fill=orange, tikzit shape=rectangle]
\tikzstyle{astate}=[fill=white, draw=black, shape=NEtriangle, tikzit category=asymmetric, tikzit shape=circle, tikzit fill=yellow]
\tikzstyle{astate conj}=[fill=white, draw=black, shape=NWtriangle, tikzit category=asymmetric, tikzit shape=circle, tikzit fill=green]
\tikzstyle{astate adj}=[fill=white, draw=black, shape=SEtriangle, tikzit category=asymmetric, tikzit shape=circle, tikzit fill=red]
\tikzstyle{astate trans}=[fill=white, draw=black, shape=SWtriangle, tikzit category=asymmetric, tikzit shape=circle, tikzit fill=orange]

\tikzstyle{hadamard edge}=[-, dashed, dash pattern=on 2pt off 0.5pt, thick, draw={rgb,255: red,68; green,136; blue,255}]
\tikzstyle{star edge}=[-, dashed, dash pattern=on 2pt off 0.5pt, thick, draw={rgb,255: red,255; green,136; blue,68}]
\tikzstyle{box edge}=[-, dashed, dash pattern=on 2pt off 0.5pt, thick, draw={rgb,255: red,203; green,192; blue,225}]
\tikzstyle{brace edge}=[-, tikzit draw=blue, decorate, decoration={brace,amplitude=1mm,raise=-1mm}]
\tikzstyle{diredge}=[->]
\tikzstyle{double edge}=[-, double, shorten <=-1mm, shorten >=-1mm, double distance=2pt]
\tikzstyle{gray edge}=[-, {gray!60!white}]
\tikzstyle{pointer edge}=[->, very thick, gray]
\tikzstyle{boldedge}=[-, line width=1.6pt, shorten <=-0.17mm, shorten >=-0.17mm]
\tikzstyle{bidir edge}=[<->, very thick, draw={rgb,255: red,191; green,191; blue,191}]

\begin{document}
\maketitle

\begin{abstract}
DisCoPy is a Python toolkit for computing with monoidal categories.
It comes with two flexible data structures for string diagrams: the first one for planar monoidal categories based on lists of layers, the second one for symmetric monoidal categories based on cospans of hypergraphs.
Algorithms for functor application then allow to translate string diagrams into code for numerical computation, be it differentiable, probabilistic or quantum.
This report gives an overview of the library and the new developments released in its version 1.0.
In particular, we showcase the implementation of diagram equality for a large fragment of the hierarchy of graphical languages for monoidal categories, as well as a new syntax for defining string diagrams as Python functions.
\end{abstract}

\section*{Extended Abstract}

String diagrams are an intuitive yet formal graphical language which has been reinvented multiple times in the context of philosophical logic~\cite{Peirce06}, circuit design~\cite{Hotz65} and spin networks~\cite{Penrose71}.
More recently, string diagrams have known a new wave of applications including quantum computing~\cite{Coecke05}, linguistics~\cite{ClarkEtAl08}, Bayesian inference~\cite{CoeckeSpekkens12}, chemical reaction networks~\cite{BaezPollard17}, databases~\cite{BonchiEtAl18}, game theory~\cite{GhaniEtAl18} and machine learning~\cite{FongEtAl19}.
Created in order to foster the development of such applications, DisCoPy is a software package that provides 1) \emph{string diagrams as a data structure} together with algorithms for composing, rewriting, drawing and checking equality between them, 2) \emph{monoidal functors for evaluating string diagrams as code}, be it for a quantum circuit, a probabilistic program or a neural network.

DisCoPy is free software, it comes with an extensive documentation and demonstration notebooks.\footnote{\url{https://docs.discopy.org}}
The library is already the topic of two tool papers~\cite{DeFeliceEtAl20,ToumiEtAl22} aimed at applied category theorists and quantum computer scientists, respectively.
It is also documented by the PhD theses of the last and first authors of this report.
The former~\cite{DeFelice22} develops a category-theoretic framework for natural language processing while the latter~\cite{Toumi22} applies this framework to \emph{quantum natural language processing}, an application which has now grown into its own library: lambeq~\cite{KartsaklisEtAl21}.
More recently, DisCoPyro~\cite{SenneshEtAl23} applied our toolkit to probabilistic generative modeling.

DisCoPy aims at becoming the fundamental package for all the applications of string diagrams.
The use of Python for applied category theory is motivated by two main reasons.
First, Python has become the programming language of reference for machine learning and quantum computing, two killer applications of category theory.
Second, Python is a programming language of choice for students and beginners.
We believe that DisCoPy can help both applied category theorists pick up programming skills and Python programmers pick up category theory concepts.
In particular, the library makes extensive use of the \emph{factory method pattern}~\cite[p.~87]{GammaEtAl95} which allows users to easily define their own custom categories.

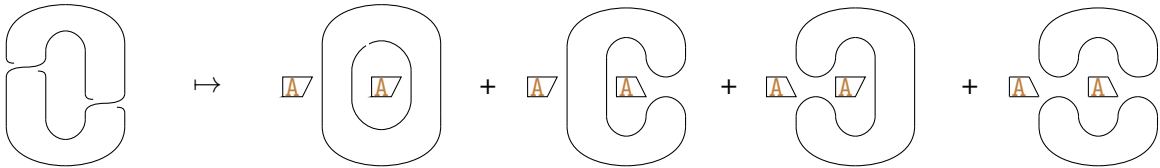
\begin{figure}[H]
\caption{Checking the equality of two symmetric diagrams by converting them to hypergraphs.}
\label{fig:equality}
\centering
\begin{minipage}{0.57\textwidth}
\begin{minted}{python}
from discopy.symmetric import Ty, Box, Swap, Diagram

x, y, z = Ty('x'), Ty('y'), Ty('z')

f = Box('f', x, y @ z)
g, h = Box('g', z, x), Box('h', y, z)

diagram_left = f >> Swap(y, z) >> g @ h
diagram_right = f >> h @ g >> Swap(z, x)

assert diagram_left != diagram_right

with Diagram.hypergraph_equality:
    assert diagram_left == diagram_right
\end{minted}
\end{minipage}
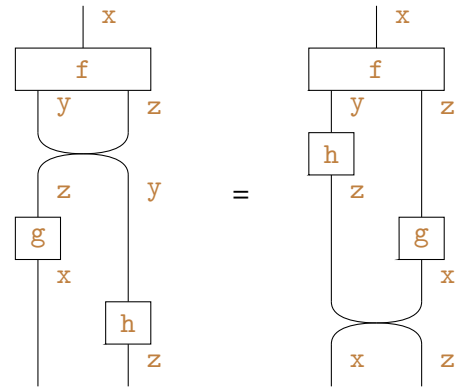
\begin{minipage}{0.38\textwidth}
\[\begin{tikzpicture}[xscale=1.6, yscale=1.5]
	\begin{pgfonlayer}{nodelayer}
		\node [style=none] (0) at (0.25, 2.25) {};
		\node [style=none] (1) at (0.75, 4.5) {};
		\node [style=none] (2) at (0.75, 4) {};
		\node [style=none, fill=white, right] (3) at (0.95, 4.375) {\py{x}};
		\node [style=none] (4) at (0.25, 3.5) {};
		\node [style=none] (5) at (0.25, 3) {};
		\node [style=none, fill=white, right] (6) at (0.45, 3.3) {\py{y}};
		\node [style=none] (7) at (1.25, 3.5) {};
		\node [style=none] (8) at (1.25, 3) {};
		\node [style=none, fill=white, right] (9) at (1.45, 3.3) {\py{z}};
		\node [style=none] (10) at (0.75, 2.75) {};
		\node [style=none] (11) at (0.25, 2.5) {};
		\node [style=none] (12) at (1.25, 2.5) {};
		\node [style=none] (13) at (0.25, 2) {};
		\node [style=none, fill=white, right] (14) at (0.45, 2.3) {\py{z}};
		\node [style=none] (15) at (1.25, 1) {};
		\node [style=none, fill=white, right] (16) at (1.45, 2.3) {\py{y}};
		\node [style=none] (17) at (0.25, 1.5) {};
		\node [style=none] (18) at (0.25, 0) {};
		\node [style=none, fill=white, right] (19) at (0.45, 1.3) {\py{x}};
		\node [style=none] (20) at (1.25, 0.5) {};
		\node [style=none] (21) at (1.25, 0) {};
		\node [style=none, fill=white, right] (22) at (1.45, 0.3) {\py{z}};
		\node [style=none] (23) at (0, 3.5) {};
		\node [style=none] (24) at (1.5, 3.5) {};
		\node [style=none] (25) at (1.5, 4) {};
		\node [style=none] (26) at (0, 4) {};
		\node [style=none, fill=white] (27) at (0.75, 3.75) {\py{f}};
		\node [style=none] (28) at (0, 1.5) {};
		\node [style=none] (29) at (0.5, 1.5) {};
		\node [style=none] (30) at (0.5, 2) {};
		\node [style=none] (31) at (0, 2) {};
		\node [style=none, fill=white] (32) at (0.25, 1.75) {\py{g}};
		\node [style=none] (33) at (1, 0.5) {};
		\node [style=none] (34) at (1.5, 0.5) {};
		\node [style=none] (35) at (1.5, 1) {};
		\node [style=none] (36) at (1, 1) {};
		\node [style=none, fill=white] (37) at (1.25, 0.75) {\py{h}};
		\node [style=none, fill=white] (38) at (2.5, 2.25) {=};
		\node [style=none] (39) at (4, 4.5) {};
		\node [style=none] (40) at (4, 4) {};
		\node [style=none, fill=white, right] (41) at (4.2, 4.375) {\py{x}};
		\node [style=none] (42) at (3.5, 3.5) {};
		\node [style=none] (43) at (3.5, 3) {};
		\node [style=none, fill=white, right] (44) at (3.7, 3.3) {\py{y}};
		\node [style=none] (45) at (4.5, 3.5) {};
		\node [style=none] (46) at (4.5, 2) {};
		\node [style=none, fill=white, right] (47) at (4.7, 3.3) {\py{z}};
		\node [style=none] (48) at (3.5, 2.5) {};
		\node [style=none] (49) at (3.5, 1) {};
		\node [style=none, fill=white, right] (50) at (3.7, 2.3) {\py{z}};
		\node [style=none] (51) at (4.5, 1.5) {};
		\node [style=none] (52) at (4.5, 1) {};
		\node [style=none, fill=white, right] (53) at (4.7, 1.3) {\py{x}};
		\node [style=none] (54) at (4, 0.75) {};
		\node [style=none] (55) at (3.5, 0.5) {};
		\node [style=none] (56) at (4.5, 0.5) {};
		\node [style=none] (57) at (3.5, 0) {};
		\node [style=none, fill=white, right] (58) at (3.7, 0.3) {\py{x}};
		\node [style=none] (59) at (4.5, 0) {};
		\node [style=none, fill=white, right] (60) at (4.7, 0.3) {\py{z}};
		\node [style=none] (61) at (3.25, 3.5) {};
		\node [style=none] (62) at (4.75, 3.5) {};
		\node [style=none] (63) at (4.75, 4) {};
		\node [style=none] (64) at (3.25, 4) {};
		\node [style=none, fill=white] (65) at (4, 3.75) {\py{f}};
		\node [style=none] (66) at (3.25, 2.5) {};
		\node [style=none] (67) at (3.75, 2.5) {};
		\node [style=none] (68) at (3.75, 3) {};
		\node [style=none] (69) at (3.25, 3) {};
		\node [style=none, fill=white] (70) at (3.5, 2.75) {\py{h}};
		\node [style=none] (71) at (4.25, 1.5) {};
		\node [style=none] (72) at (4.75, 1.5) {};
		\node [style=none] (73) at (4.75, 2) {};
		\node [style=none] (74) at (4.25, 2) {};
		\node [style=none, fill=white] (75) at (4.5, 1.75) {\py{g}};
	\end{pgfonlayer}
	\begin{pgfonlayer}{edgelayer}
		\draw [in=90, out=-90] (1.center) to (2.center);
		\draw [in=90, out=-90] (4.center) to (5.center);
		\draw [in=90, out=-90] (7.center) to (8.center);
		\draw [in=90, out=180, looseness=0.94] (10.center) to (11.center);
		\draw [in=90, out=0, looseness=0.94] (10.center) to (12.center);
		\draw [in=180, out=-90, looseness=0.94] (5.center) to (10.center);
		\draw [in=0, out=-90, looseness=0.94] (8.center) to (10.center);
		\draw [in=90, out=-90] (11.center) to (13.center);
		\draw [in=90, out=-90] (12.center) to (15.center);
		\draw [in=90, out=-90] (17.center) to (18.center);
		\draw [in=90, out=-90] (20.center) to (21.center);
		\draw [-, fill=white] (23.center)
			 to (24.center)
			 to (25.center)
			 to (26.center)
			 to cycle;
		\draw [-, fill=white] (28.center)
			 to (29.center)
			 to (30.center)
			 to (31.center)
			 to cycle;
		\draw [-, fill=white] (33.center)
			 to (34.center)
			 to (35.center)
			 to (36.center)
			 to cycle;
		\draw [in=90, out=-90] (39.center) to (40.center);
		\draw [in=90, out=-90] (42.center) to (43.center);
		\draw [in=90, out=-90] (45.center) to (46.center);
		\draw [in=90, out=-90] (48.center) to (49.center);
		\draw [in=90, out=-90] (51.center) to (52.center);
		\draw [in=90, out=180, looseness=0.94] (54.center) to (55.center);
		\draw [in=90, out=0, looseness=0.94] (54.center) to (56.center);
		\draw [in=180, out=-90, looseness=0.94] (49.center) to (54.center);
		\draw [in=0, out=-90, looseness=0.94] (52.center) to (54.center);
		\draw [in=90, out=-90] (55.center) to (57.center);
		\draw [in=90, out=-90] (56.center) to (59.center);
		\draw [-, fill=white] (61.center)
			 to (62.center)
			 to (63.center)
			 to (64.center)
			 to cycle;
		\draw [-, fill=white] (66.center)
			 to (67.center)
			 to (68.center)
			 to (69.center)
			 to cycle;
		\draw [-, fill=white] (71.center)
			 to (72.center)
			 to (73.center)
			 to (74.center)
			 to cycle;
	\end{pgfonlayer}
\end{tikzpicture}\]
\end{minipage}
\end{figure}

So what is DisCoPy?
In a nutshell, it is a \emph{domain specific language} (DSL) for morphisms in (pre)monoidal categories.
Its main data structure is \py{Diagram}, an implementation of the arrows of the free premonoidal category generated by the class \py{Ob} as objects and the class \py{Box} as arrows.
These can be constructed either using a point-free syntax with composition (\py{>>}) and tensor (\py{@}) as in Figure~\ref{fig:equality} or with the standard syntax for Python functions as demonstrated in Figure~\ref{fig:copy-discard}.
Its main algorithm is \py{Functor} application which evaluates diagrams as morphisms in an arbitrary (pre)monoidal category.
Endofunctors on the free premonoidal category allow to ``open the box'' by replacing it with an arbitrary diagram.

With \py{python}, \py{matrix} or \py{tensor} as codomain, functors allow to turn diagrams into fast numerical computation via interfaces with any of NumPy~\cite{VanDerWaltEtAl11}, TensorFlow~\cite{AbadiEtAl16}, PyTorch~\cite{PaszkeEtAl19}, JAX~\cite{BradburyEtAl18} or TensorNetwork~\cite{RobertsEtAl19}.
The \py{grammar} subpackage interfaces with parsers for natural language processing while the \py{quantum} subpackage interfaces with tket~\cite{SivarajahEtAl20} for circuit compilation and PyZX~\cite{KissingerVanDeWetering19} for rewriting.
It also implements the classical simulation of quantum circuits as unitary matrices or quantum channels as well as \emph{diagrammatic differentiation}~\cite{ToumiEtAl21}, i.e. automatic differentiation of parameterised string diagrams.

With the release of its version 1.0, the library has undergone a complete refactor which simplifies its architecture while making it more modular.
In particular, \py{Diagram} is now a subclass of \py{Arrow} (the implementation of the free category) so that identity and composition are defined exactly once.
Another important change is that \py{Matrix} and \py{Tensor} now are generic types parameterised by their data type.

Among the new features of this v1.0, the most significant is the implementation of a large fragment of the hierarchy of graphical languages for monoidal categories as surveyed by Selinger~\cite{Selinger10}, which is summarised in Figure~\ref{hierarchy}.
Each module implements a DSL for morphisms in monoidal categories with extra structure (e.g. \py{braided}, \py{traced}, \py{closed}, etc.) with its own subclasses of \py{Diagram} and \py{Functor}.
In the cases when it has a known solution, we have implemented the \emph{word problem} for the free categories, i.e. decide whether two diagrams are equal up to the axiom of the category.

In the case of symmetric, traced, compact and hypergraph categories, this reduces to \emph{hypergraph isomorphism} which we compute via the graph isomorphism algorithm of NetworkX~\cite{HagbergEtAl08}.
This is implemented in a new \py{Hypergraph} class which provides an alternative representation of string diagrams where the axioms for symmetric categories and special commutative Frobenius algebras hold for free.

Plans for future developments include the implementation of free bicategories in terms of diagrams with colours, as well as double categories where wires can go horizontally.
Another promising direction is the implementation of \emph{double-pushout rewriting} via interfaces to existing libraries~\cite{AndersenEtAl16,SobocinskiEtAl19,HarmerOshurko20}.
Diagram rewriting can then itself be represented in terms of \emph{higher-dimensional diagrams}, with a generalisation of our list-based \py{Diagram} data structure~\cite{BarVicary17} or with \emph{combinatorial directed cell complexes}~\cite{HadzihasanovicKessler22} which generalise our graph-based \py{Hypergraph} data structure.

The appendix of this report is structured as follows: Section~\ref{cat} covers the basics of free categories and functors, Section~\ref{monoidal} and \ref{hypergraph} introduce our two data structures for planar diagrams and hypergraphs.
\vfill
\pagebreak

\begin{figure}[H]
\caption{Architecture of the DisCoPy library.}
\vspace{-10pt}
\label{hierarchy}
\[\begin{tikzpicture}[scale=1.33, every node/.style={scale=1.33}]
	\begin{pgfonlayer}{nodelayer}
		\node [style=none] (0) at (0, 2) {\pyref{cat}{cat}};
		\node [style=none] (1) at (0, 0) {\pyref{monoidal}{monoidal}$\ssymbol{1}$};
		\node [style=none] (2) at (0, -2) {\pyref{traced}{traced}$\ssymbol{2}$};
		\node [style=none] (3) at (-4, -2) {\pyref{braided}{braided}$\ssymbol{3}$};
		\node [style=none] (4) at (4, -2) {\pyref{closed}{closed}};
		\node [style=none] (5) at (4, -4) {\pyref{rigid}{rigid}$\ssymbol{3}$};
		\node [style=none] (6) at (4, -6) {\pyref{pivotal}{pivotal}};
		\node [style=none] (7) at (4, -8) {\pyref{ribbon}{ribbon}$\ssymbol{4}$};
		\node [style=none] (8) at (-4, -8) {\pyref{balanced}{balanced}$\ssymbol{4}$};
		\node [style=none] (9) at (-4, -10) {\pyref{symmetric}{symmetric}$\ssymbol{5}$};
		\node [style=none] (10) at (4, -10) {\pyref{compact}{compact}$\ssymbol{5}$};
		\node [style=none] (11) at (-4, -12) {\pyref{comonoid}{comonoid}$\ssymbol{5}$};
		\node [style=none] (12) at (0, -11) {\pyref{hypergraph}{hypergraph}};
		\node [style=none] (13) at (4, -12) {\pyref{frobenius}{frobenius}$\ssymbol{5}$};
		\node [style=none] (14) at (0, -0.5) {};
		\node [style=none] (15) at (0, -1.5) {};
		\node [style=none] (16) at (0, 1.5) {};
		\node [style=none] (17) at (0, 0.5) {};
		\node [style=none] (18) at (1, -0.5) {};
		\node [style=none] (19) at (3, -1.5) {};
		\node [style=none] (20) at (-1, -0.5) {};
		\node [style=none] (21) at (-3, -1.5) {};
		\node [style=none] (22) at (-4, -2.5) {};
		\node [style=none] (23) at (-4, -7.5) {};
		\node [style=none] (24) at (4, -2.5) {};
		\node [style=none] (25) at (4, -3.5) {};
		\node [style=none] (26) at (4, -4.5) {};
		\node [style=none] (27) at (4, -5.5) {};
		\node [style=none] (28) at (4, -6.5) {};
		\node [style=none] (29) at (4, -7.5) {};
		\node [style=none] (30) at (-4, -8.5) {};
		\node [style=none] (31) at (-4, -9.5) {};
		\node [style=none] (32) at (-4, -10.5) {};
		\node [style=none] (33) at (-4, -11.5) {};
		\node [style=none] (34) at (4, -10.5) {};
		\node [style=none] (35) at (4, -11.5) {};
		\node [style=none] (36) at (4, -8.5) {};
		\node [style=none] (37) at (4, -9.5) {};
		\node [style=none] (38) at (-2.5, -10) {};
		\node [style=none] (39) at (2.5, -10) {};
		\node [style=none] (40) at (-2.5, -8) {};
		\node [style=none] (41) at (2.5, -8) {};
		\node [style=none] (44) at (3.5, -5.5) {};
		\node [style=none] (45) at (0.5, -2.75) {};
		\node [style=none] (46) at (-2.5, -12) {};
		\node [style=none] (47) at (2.5, -12) {};
		\node [style=none] (48) at (0, -9) {\pyref{interaction}{interaction}};
		\node [style=none] (49) at (-4, -15) {\pyref{python}{python}};
		\node [style=none] (50) at (0, -15) {\pyref{matrix}{matrix}};
		\node [style=none] (51) at (4, -15) {\pyref{tensor}{tensor}};
		\node [style=none] (52) at (9, 0) {\pyref{cfg}{cfg}};
		\node [style=none] (54) at (9, -2) {\pyref{categorial}{categorial}};
		\node [style=none] (55) at (9, -4) {\pyref{pregroup}{pregroup}};
		\node [style=none] (64) at (9, -6) {\pyref{dependency}{dependency}};
		\node [style=none] (65) at (9, -9) {\pyref{circuit}{circuit}};
		\node [style=none] (66) at (9, -11) {\pyref{tk}{tk}};
		\node [style=none] (67) at (9, -13) {\pyref{zx}{zx}};
		\node [style=none] (68) at (9, -15) {\pyref{channel}{channel}};
		\node [style=none] (78) at (-5.75, -14) {};
		\node [style=none] (79) at (6, -14.25) {};
		\node [style=none] (80) at (6, -15.75) {};
		\node [style=none] (81) at (-6, -15.75) {};
		\node [style=none] (86) at (0, 3.375) {\href{https://docs.discopy.org/en/main/api/syntax.html}{\scriptsize\py{syntax}}};
		\node [style=none] (87) at (9, 1.375) {\href{https://docs.discopy.org/en/main/api/grammar.html}{\scriptsize\py{grammar}}};
		\node [style=none] (88) at (9, -7.625) {\href{https://docs.discopy.org/en/main/api/quantum.html}{\scriptsize\py{quantum}}};
		\node [style=none] (89) at (0, -13.625) {\href{https://docs.discopy.org/en/main/api/semantics.html}{\scriptsize\py{semantics}}};
		\node [style=none] (90) at (-6, -14.25) {};
		\node [style=none] (91) at (-5.75, -16) {};
		\node [style=none] (92) at (5.75, -16) {};
		\node [style=none] (93) at (5.75, -14) {};
		\node [style=none] (94) at (-5.75, 3) {};
		\node [style=none] (95) at (6, 2.75) {};
		\node [style=none] (96) at (6, -12.75) {};
		\node [style=none] (97) at (-6, -12.75) {};
		\node [style=none] (98) at (-6, 2.75) {};
		\node [style=none] (99) at (-5.75, -13) {};
		\node [style=none] (100) at (5.75, -13) {};
		\node [style=none] (101) at (5.75, 3) {};
		\node [style=none] (102) at (7.25, 1) {};
		\node [style=none] (103) at (11, 0.75) {};
		\node [style=none] (104) at (11, -6.75) {};
		\node [style=none] (105) at (7, -6.75) {};
		\node [style=none] (106) at (7, 0.75) {};
		\node [style=none] (107) at (7.25, -7) {};
		\node [style=none] (108) at (10.75, -7) {};
		\node [style=none] (109) at (10.75, 1) {};
		\node [style=none] (110) at (7.25, -8) {};
		\node [style=none] (111) at (11, -8.25) {};
		\node [style=none] (112) at (11, -15.75) {};
		\node [style=none] (113) at (7, -15.75) {};
		\node [style=none] (114) at (7, -8.25) {};
		\node [style=none] (115) at (7.25, -16) {};
		\node [style=none] (116) at (10.75, -16) {};
		\node [style=none] (117) at (10.75, -8) {};
		\node [style=none] (118) at (-0.5, -2.75) {};
		\node [style=none] (119) at (-3.5, -7.5) {};
	\end{pgfonlayer}
	\begin{pgfonlayer}{edgelayer}
		\draw [style=diredge] (16.center) to (17.center);
		\draw [style=diredge] (14.center) to (15.center);
		\draw [style=diredge] (18.center) to (19.center);
		\draw [style=diredge] (20.center) to (21.center);
		\draw [style=diredge] (22.center) to (23.center);
		\draw [style=diredge] (24.center) to (25.center);
		\draw [style=diredge] (26.center) to (27.center);
		\draw [style=diredge] (28.center) to (29.center);
		\draw [style=diredge] (30.center) to (31.center);
		\draw [style=diredge] (32.center) to (33.center);
		\draw [style=diredge] (34.center) to (35.center);
		\draw [style=diredge] (36.center) to (37.center);
		\draw [style=diredge] (38.center) to (39.center);
		\draw [style=diredge] (40.center) to (41.center);
		\draw [style=diredge] (45.center) to (44.center);
		\draw [style=diredge] (46.center) to (47.center);
		\draw [style=diredge] (118.center) to (119.center);
		\draw (80.center)
			 to [bend left=45, looseness=0.75] (92.center)
			 to (91.center)
			 to [bend left=45, looseness=0.75] (81.center)
			 to (90.center)
			 to [bend left=45, looseness=0.75] (78.center)
			 to (93.center)
			 to [bend left=45, looseness=0.75] (79.center)
			 to cycle;
		\draw (96.center)
			 to [bend left=45, looseness=0.75] (100.center)
			 to (99.center)
			 to [bend left=45, looseness=0.75] (97.center)
			 to (98.center)
			 to [bend left=45, looseness=0.75] (94.center)
			 to (101.center)
			 to [bend left=45, looseness=0.75] (95.center)
			 to cycle;
		\draw (104.center)
			 to [bend left=45, looseness=0.75] (108.center)
			 to (107.center)
			 to [bend left=45, looseness=0.75] (105.center)
			 to (106.center)
			 to [bend left=45, looseness=0.75] (102.center)
			 to (109.center)
			 to [bend left=45, looseness=0.75] (103.center)
			 to cycle;
		\draw (112.center)
			 to [bend left=45, looseness=0.75] (116.center)
			 to (115.center)
			 to [bend left=45, looseness=0.75] (113.center)
			 to (114.center)
			 to [bend left=45, looseness=0.75] (110.center)
			 to (117.center)
			 to [bend left=45, looseness=0.75] (111.center)
			 to cycle;
	\end{pgfonlayer}
\end{tikzpicture}\]
\caption*
{$\ssymbol{1}$ Planar diagram equality can be computed in quadratic time and normal form in cubic time~\cite{DelpeuchVicary22}.\\
$\ssymbol{2}$ The \py{traced} module implements planar traced categories while \py{balanced} is traced by default.\\
$\ssymbol{3}$ Rigid diagram normal form is computed in polynomial time with \emph{snake removal}, see~\cite[§1.4.1]{Toumi22}.\\
$\ssymbol{4}$ Braided diagram equality is unknot hard,~\cite{DelpeuchVicary21} it is currently not implemented in DisCoPy.\\
$\ssymbol{5}$ Hypergraph diagram equality reduces to graph isomorphism, it is implemented with NetworkX~\cite{HagbergEtAl08}.}
\end{figure}

\begin{figure}[H]
\caption{Defining a diagram as a list of layers.}
\label{fig:basic-diagram}
\centering
\begin{minipage}{0.65\textwidth}
\begin{minted}{python}
from discopy.monoidal import Ty, Box, Layer, Diagram

x, y, z = Ty('x'), Ty('y'), Ty('z')
f, g, h = Box('f', x, y @ z), Box('g', y, z), Box('h', z, z)

assert f >> g @ h == Diagram(
    dom=x, cod=z @ z, inside=(
        Layer(Ty(), f, Ty()),
        Layer(Ty(), g, z),
        Layer(z,    h, Ty())))
\end{minted}
\end{minipage}
\begin{minipage}[b]{0.33\textwidth}
\[\begin{tikzpicture}[baseline={(0.base)}]
	\begin{pgfonlayer}{nodelayer}
		\node [style=none] (0) at (0.25, 2.625) {};
		\node [style=none] (1) at (1, 5.375) {};
		\node [style=none] (2) at (1, 4.625) {};
		\node [style=none, fill=white, right] (3) at (1.1, 5.125) {\py{x}};
		\node [style=none] (4) at (0.25, 3.875) {};
		\node [style=none] (5) at (0.25, 3) {};
		\node [style=none, fill=white, right] (6) at (0.35, 3.525) {\py{y}};
		\node [style=none] (7) at (1.75, 3.875) {};
		\node [style=none] (8) at (1.75, 1.625) {};
		\node [style=none, fill=white, right] (9) at (1.85, 3.525) {\py{z}};
		\node [style=none] (10) at (0.25, 2.25) {};
		\node [style=none] (11) at (0.25, 0) {};
		\node [style=none, fill=white, right] (12) at (0.35, 1.9) {\py{z}};
		\node [style=none] (13) at (1.75, 0.875) {};
		\node [style=none] (14) at (1.75, 0) {};
		\node [style=none, fill=white, right] (15) at (1.85, 0.525) {\py{z}};
		\node [style=none] (16) at (-0.125, 3.875) {};
		\node [style=none] (17) at (2.125, 3.875) {};
		\node [style=none] (18) at (2.125, 4.625) {};
		\node [style=none] (19) at (-0.125, 4.625) {};
		\node [style=none, fill=white] (20) at (1, 4.25) {\py{f}};
		\node [style=none] (21) at (-0.125, 2.25) {};
		\node [style=none] (22) at (0.625, 2.25) {};
		\node [style=none] (23) at (0.625, 3) {};
		\node [style=none] (24) at (-0.125, 3) {};
		\node [style=none, fill=white] (25) at (0.25, 2.625) {\py{g}};
		\node [style=none] (26) at (1.375, 0.875) {};
		\node [style=none] (27) at (2.125, 0.875) {};
		\node [style=none] (28) at (2.125, 1.625) {};
		\node [style=none] (29) at (1.375, 1.625) {};
		\node [style=none, fill=white] (30) at (1.75, 1.25) {\py{h}};
	\end{pgfonlayer}
	\begin{pgfonlayer}{edgelayer}
		\draw [in=90, out=-90] (1.center) to (2.center);
		\draw [in=90, out=-90] (4.center) to (5.center);
		\draw [in=90, out=-90] (7.center) to (8.center);
		\draw [in=90, out=-90] (10.center) to (11.center);
		\draw [in=90, out=-90] (13.center) to (14.center);
		\draw [-, fill=white] (16.center)
			 to (17.center)
			 to (18.center)
			 to (19.center)
			 to cycle;
		\draw [-, fill=white] (21.center)
			 to (22.center)
			 to (23.center)
			 to (24.center)
			 to cycle;
		\draw [-, fill=white] (26.center)
			 to (27.center)
			 to (28.center)
			 to (29.center)
			 to cycle;
	\end{pgfonlayer}
\end{tikzpicture}\]
\end{minipage}
\end{figure}

\begin{figure}[H]
\caption{Defining Boolean circuits as a subclass of \py{Diagram} with natural numbers as objects.}
\label{fig:circuits}
\centering
\begin{minipage}{0.73\textwidth}
\begin{minted}{python}
from discopy import monoidal, python
from discopy.cat import factory, Category

@factory  # Ensure that composition of circuits remains a circuit.
class Circuit(monoidal.Diagram):
    ty_factory = monoidal.PRO  # Use natural numbers as objects.

    def __call__(self, *bits):
        F = monoidal.Functor(
            ob=lambda _: (bool, ), ar=lambda f: f.data,
            cod=Category(python.Ty, python.Function))
        return F(self)(bits)

class Gate(monoidal.Box, Circuit):
    """A gate is just a box in a circuit with a function as data."""

NAND = Gate("NAND", 2, 1, data=lambda x, y: not (x and y))
COPY = Gate("COPY", 1, 2, data=lambda x: (x, x))

XOR = COPY @ COPY >> 1 @ (NAND >> COPY) @ 1 >> NAND @ NAND >> NAND
CNOT = COPY @ 1 >> 1 @ XOR
NOTC = 1 @ COPY >> XOR @ 1
SWAP = CNOT >> NOTC >> CNOT  # Exercise: Find a cheaper SWAP circuit!

assert all(SWAP(x, y) == (y, x) for x in [True, False]
                                for y in [True, False])
\end{minted}
\end{minipage}
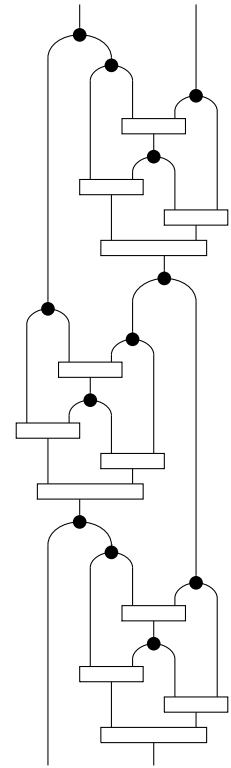
\begin{minipage}{0.26\textwidth}
\vfill \[\input{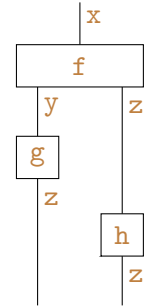}\] \vfill
\end{minipage}
\end{figure}

\begin{figure}
\caption{Spiral-shaped diagrams are the worst-case for normalising planar string diagrams.}
\label{spiral-example}
\centering
\begin{minted}{python}
from discopy.monoidal import Ty, Box
from discopy.drawing import Equation

x = Ty('x')
f, u = Box('f', Ty(), x @ x), Box('u', Ty(), x)

def spiral(length):
    diagram, n = u, length // 2 - 1
    for i in range(n):
        diagram >>= x ** i @ f @ x ** (i + 1)
    diagram >>= x ** n @ u.dagger() @ x ** n
    for i in range(n):
        m = n - i - 1
        diagram >>= x ** m @ f.dagger() @ x ** m
    return diagram

assert spiral(8).dagger() != spiral(8)
assert spiral(8).dagger() == spiral(8).normal_form()
\end{minted}

\[\begin{tikzpicture}[baseline=(0.base), yscale=0.7, xscale=.8]
\begin{pgfonlayer}{nodelayer}
\node (0) at (0, 4.0) {};
\node [] (1) at (6.0, 6.75) {};
\node [] (2) at (6.0, 0.25) {};
\node [] (3) at (0.0, 5.75) {};
\node [] (4) at (0.0, 0.25) {};
\node [] (5) at (5.0, 5.75) {};
\node [] (6) at (5.0, 1.25) {};
\node [] (7) at (1.0, 4.75) {};
\node [] (8) at (1.0, 1.25) {};
\node [] (9) at (4.0, 4.75) {};
\node [] (10) at (4.0, 2.25) {};
\node [] (11) at (2.0, 3.75) {};
\node [] (12) at (2.0, 2.25) {};
\node [] (13) at (3.0, 3.75) {};
\node [] (14) at (3.0, 3.25) {};
\node [] (15) at (5.75, 6.75) {};
\node [] (16) at (6.25, 6.75) {};
\node [] (17) at (6.5, 7.25) {};
\node [] (18) at (5.75, 7.25) {};
\node [style=none, fill=white] (19) at (6.0, 7.0) {};
\node [] (20) at (-0.25, 5.75) {};
\node [] (21) at (5.25, 5.75) {};
\node [] (22) at (5.5, 6.25) {};
\node [] (23) at (-0.25, 6.25) {};
\node [style=none, fill=white] (24) at (2.5, 6.0) {};
\node [] (25) at (0.75, 4.75) {};
\node [] (26) at (4.25, 4.75) {};
\node [] (27) at (4.5, 5.25) {};
\node [] (28) at (0.75, 5.25) {};
\node [style=none, fill=white] (29) at (2.5, 5.0) {};
\node [] (30) at (1.75, 3.75) {};
\node [] (31) at (3.25, 3.75) {};
\node [] (32) at (3.5, 4.25) {};
\node [] (33) at (1.75, 4.25) {};
\node [style=none, fill=white] (34) at (2.5, 4.0) {};
\node [] (35) at (2.75, 2.75) {};
\node [] (36) at (3.5, 2.75) {};
\node [] (37) at (3.25, 3.25) {};
\node [] (38) at (2.75, 3.25) {};
\node [style=none, fill=white] (39) at (3.0, 3.0) {};
\node [] (40) at (1.75, 1.75) {};
\node [] (41) at (4.5, 1.75) {};
\node [] (42) at (4.25, 2.25) {};
\node [] (43) at (1.75, 2.25) {};
\node [style=none, fill=white] (44) at (3.0, 2.0) {};
\node [] (45) at (0.75, 0.75) {};
\node [] (46) at (5.5, 0.75) {};
\node [] (47) at (5.25, 1.25) {};
\node [] (48) at (0.75, 1.25) {};
\node [style=none, fill=white] (49) at (3.0, 1.0) {};
\node [] (50) at (-0.25, -0.25) {};
\node [] (51) at (6.5, -0.25) {};
\node [] (52) at (6.25, 0.25) {};
\node [] (53) at (-0.25, 0.25) {};
\node [style=none, fill=white] (54) at (3.0, 0.0) {};
\node [style=none, fill=white] (55) at (7.5, 4.0) {=};
\node [] (56) at (8.5, 6.75) {};
\node [] (57) at (8.5, 1.25) {};
\node [] (58) at (14.5, 6.75) {};
\node [] (59) at (14.5, 0.25) {};
\node [] (60) at (9.5, 5.75) {};
\node [] (61) at (9.5, 2.25) {};
\node [] (62) at (13.5, 5.75) {};
\node [] (63) at (13.5, 1.25) {};
\node [] (64) at (10.5, 4.75) {};
\node [] (65) at (10.5, 3.25) {};
\node [] (66) at (12.5, 4.75) {};
\node [] (67) at (12.5, 2.25) {};
\node [] (68) at (11.5, 3.75) {};
\node [] (69) at (11.5, 3.25) {};
\node [] (70) at (8.25, 6.75) {};
\node [] (71) at (14.75, 6.75) {};
\node [] (72) at (15.0, 7.25) {};
\node [] (73) at (8.25, 7.25) {};
\node [style=none, fill=white] (74) at (11.5, 7.0) {};
\node [] (75) at (9.25, 5.75) {};
\node [] (76) at (13.75, 5.75) {};
\node [] (77) at (14.0, 6.25) {};
\node [] (78) at (9.25, 6.25) {};
\node [style=none, fill=white] (79) at (11.5, 6.0) {};
\node [] (80) at (10.25, 4.75) {};
\node [] (81) at (12.75, 4.75) {};
\node [] (82) at (13.0, 5.25) {};
\node [] (83) at (10.25, 5.25) {};
\node [style=none, fill=white] (84) at (11.5, 5.0) {};
\node [] (85) at (11.25, 3.75) {};
\node [] (86) at (11.75, 3.75) {};
\node [] (87) at (12.0, 4.25) {};
\node [] (88) at (11.25, 4.25) {};
\node [style=none, fill=white] (89) at (11.5, 4.0) {};
\node [] (90) at (10.25, 2.75) {};
\node [] (91) at (12.0, 2.75) {};
\node [] (92) at (11.75, 3.25) {};
\node [] (93) at (10.25, 3.25) {};
\node [style=none, fill=white] (94) at (11.0, 3.0) {};
\node [] (95) at (9.25, 1.75) {};
\node [] (96) at (13.0, 1.75) {};
\node [] (97) at (12.75, 2.25) {};
\node [] (98) at (9.25, 2.25) {};
\node [style=none, fill=white] (99) at (11.0, 2.0) {};
\node [] (100) at (8.25, 0.75) {};
\node [] (101) at (14.0, 0.75) {};
\node [] (102) at (13.75, 1.25) {};
\node [] (103) at (8.25, 1.25) {};
\node [style=none, fill=white] (104) at (11.0, 1.0) {};
\node [] (105) at (14.25, -0.25) {};
\node [] (106) at (15.0, -0.25) {};
\node [] (107) at (14.75, 0.25) {};
\node [] (108) at (14.25, 0.25) {};
\node [style=none, fill=white] (109) at (14.5, 0.0) {};
\end{pgfonlayer}
\begin{pgfonlayer}{edgelayer}
\draw [in=90, out=-90] (1.center) to (2.center);
\draw [in=90, out=-90] (3.center) to (4.center);
\draw [in=90, out=-90] (5.center) to (6.center);
\draw [in=90, out=-90] (7.center) to (8.center);
\draw [in=90, out=-90] (9.center) to (10.center);
\draw [in=90, out=-90] (11.center) to (12.center);
\draw [in=90, out=-90] (13.center) to (14.center);
\draw [-, fill={white}] (15.center) to (16.center) to (17.center) to (18.center) to (15.center);
\draw [-, fill={white}] (20.center) to (21.center) to (22.center) to (23.center) to (20.center);
\draw [-, fill={white}] (25.center) to (26.center) to (27.center) to (28.center) to (25.center);
\draw [-, fill={white}] (30.center) to (31.center) to (32.center) to (33.center) to (30.center);
\draw [-, fill={white}] (35.center) to (36.center) to (37.center) to (38.center) to (35.center);
\draw [-, fill={white}] (40.center) to (41.center) to (42.center) to (43.center) to (40.center);
\draw [-, fill={white}] (45.center) to (46.center) to (47.center) to (48.center) to (45.center);
\draw [-, fill={white}] (50.center) to (51.center) to (52.center) to (53.center) to (50.center);
\draw [in=90, out=-90] (56.center) to (57.center);
\draw [in=90, out=-90] (58.center) to (59.center);
\draw [in=90, out=-90] (60.center) to (61.center);
\draw [in=90, out=-90] (62.center) to (63.center);
\draw [in=90, out=-90] (64.center) to (65.center);
\draw [in=90, out=-90] (66.center) to (67.center);
\draw [in=90, out=-90] (68.center) to (69.center);
\draw [-, fill={white}] (70.center) to (71.center) to (72.center) to (73.center) to (70.center);
\draw [-, fill={white}] (75.center) to (76.center) to (77.center) to (78.center) to (75.center);
\draw [-, fill={white}] (80.center) to (81.center) to (82.center) to (83.center) to (80.center);
\draw [-, fill={white}] (85.center) to (86.center) to (87.center) to (88.center) to (85.center);
\draw [-, fill={white}] (90.center) to (91.center) to (92.center) to (93.center) to (90.center);
\draw [-, fill={white}] (95.center) to (96.center) to (97.center) to (98.center) to (95.center);
\draw [-, fill={white}] (100.center) to (101.center) to (102.center) to (103.center) to (100.center);
\draw [-, fill={white}] (105.center) to (106.center) to (107.center) to (108.center) to (105.center);
\end{pgfonlayer}
\end{tikzpicture}\]
\end{figure}

\begin{figure}
\caption{Computing the golden ratio as the trace of a string diagram interpreted as a fixed point.}
\label{traced-example}
\centering
\begin{minipage}{0.65\textwidth}
\begin{minted}{python}
from discopy.traced import Ty, Box, Category, Functor
from discopy import python

x = Ty('x')
add, div = Box('+', x @ x, x), Box('/', x @ x, x)
copy, one = Box('', x, x @ x), Box('1', Ty(), x)

phi = ((one >> copy) @ x >> x @ div >> add >> copy).trace()

# The default y=1 is the initial value for the fixed point.
F = Functor(ob={x: int},
            ar={div: lambda x, y=1: x / y,
                add: lambda x, y: x + y,
                copy: lambda x: (x, x),
                one: lambda: 1},
            cod=Category(python.Ty, python.Function))

assert F(phi)() == 0.5 * (1 + 5 ** 0.5)
\end{minted}
\end{minipage}
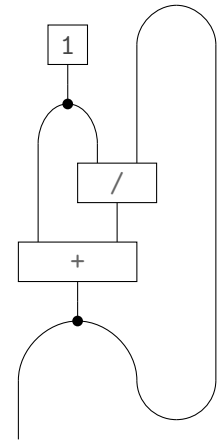
\hfill
\begin{minipage}{0.34\textwidth}
\[\begin{tikzpicture}[scale=1.4]
	\begin{pgfonlayer}{nodelayer}
		\node [style=none] (1) at (2.5, 5.75) {};
		\node [style=none] (2) at (2, 5.25) {};
		\node [style=none] (3) at (3, 5.25) {};
		\node [style=none] (4) at (2, 3.75) {};
		\node [style=none] (5) at (3, 1) {};
		\node [style=none] (6) at (1.125, 5) {};
		\node [style=none] (7) at (1.125, 4.75) {};
		\node [style=none] (8) at (1.125, 4.5) {};
		\node [style=none] (9) at (0.75, 4) {};
		\node [style=none] (10) at (1.5, 4) {};
		\node [style=none] (11) at (0.75, 2.75) {};
		\node [style=none] (12) at (1.5, 3.75) {};
		\node [style=none] (13) at (1.75, 3.25) {};
		\node [style=none] (14) at (1.75, 2.75) {};
		\node [style=none] (15) at (1.25, 2.25) {};
		\node [style=none] (16) at (1.25, 2) {};
		\node [style=none] (17) at (1.25, 1.75) {};
		\node [style=none] (18) at (0.5, 1) {};
		\node [style=none] (19) at (2, 1) {};
		\node [style=none] (20) at (0.5, 0.25) {};
		\node [circle, fill=black, scale=0.392] (23) at (1.125, 4.5) {};
		\node [circle, fill=black, scale=0.392] (23) at (1.25, 1.75) {};
		\node [style=none] (23) at (0.875, 5) {};
		\node [style=none] (24) at (1.375, 5) {};
		\node [style=none] (25) at (1.375, 5.5) {};
		\node [style=none] (26) at (0.875, 5.5) {};
		\node [style=none, fill=white] (27) at (1.125, 5.25) {\py{1}};
		\node [style=none] (28) at (1.25, 3.25) {};
		\node [style=none] (29) at (2.25, 3.25) {};
		\node [style=none] (30) at (2.25, 3.75) {};
		\node [style=none] (31) at (1.25, 3.75) {};
		\node [style=none, fill=white] (32) at (1.75, 3.5) {\py{/}};
		\node [style=none] (33) at (0.5, 2.25) {};
		\node [style=none] (34) at (2, 2.25) {};
		\node [style=none] (35) at (2, 2.75) {};
		\node [style=none] (36) at (0.5, 2.75) {};
		\node [style=none, fill=white] (37) at (1.25, 2.5) {\py{+}};
		\node [style=none] (38) at (2.5, 0.5) {};
	\end{pgfonlayer}
	\begin{pgfonlayer}{edgelayer}
		\draw [in=90, out=180, looseness=0.94] (1.center) to (2.center);
		\draw [in=90, out=0, looseness=0.94] (1.center) to (3.center);
		\draw [in=90, out=-90] (2.center) to (4.center);
		\draw [in=90, out=-90] (3.center) to (5.center);
		\draw [in=90, out=-90] (6.center) to (7.center);
		\draw [in=90, out=180, looseness=0.94] (8.center) to (9.center);
		\draw [in=90, out=0, looseness=0.94] (8.center) to (10.center);
		\draw [in=90, out=-90] (7.center) to (8.center);
		\draw [in=90, out=-90] (9.center) to (11.center);
		\draw [in=90, out=-90] (10.center) to (12.center);
		\draw [in=90, out=-90] (13.center) to (14.center);
		\draw [in=90, out=-90] (15.center) to (16.center);
		\draw [in=90, out=180, looseness=0.94] (17.center) to (18.center);
		\draw [in=90, out=0, looseness=0.94] (17.center) to (19.center);
		\draw [in=90, out=-90] (16.center) to (17.center);
		\draw [in=90, out=-90] (18.center) to (20.center);
		\draw [-, fill=white] (23.center)
			 to (24.center)
			 to (25.center)
			 to (26.center)
			 to cycle;
		\draw [-, fill=white] (28.center)
			 to (29.center)
			 to (30.center)
			 to (31.center)
			 to cycle;
		\draw [-, fill=white] (33.center)
			 to (34.center)
			 to (35.center)
			 to (36.center)
			 to cycle;
		\draw [bend right=45] (38.center) to (5.center);
		\draw [bend right=45] (19.center) to (38.center);
	\end{pgfonlayer}
\end{tikzpicture}\]
\end{minipage}
\end{figure}
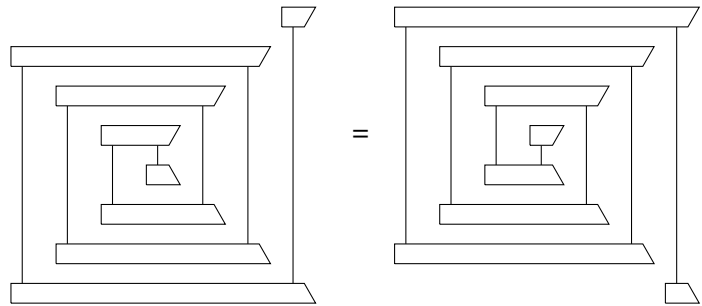

\begin{figure}[H]
\caption{The Kauffman bracket as a ribbon functor into a category where the braiding is a formal sum.}
\label{fig:kauffman}
\centering
\begin{minted}{python}
from discopy import ribbon, drawing
from discopy.cat import factory, Category

x = ribbon.Ty('x')
cup, cap, braid = ribbon.Cup(x.r, x), ribbon.Cap(x.r, x), ribbon.Braid(x, x)
link = cap >> x.r @ cap @ x >> braid.r @ braid >> x.r @ cup @ x >> cup

@factory
class Kauffman(ribbon.Diagram):
    ty_factory = ribbon.PRO

class Cup(ribbon.Cup, Kauffman): pass
class Cap(ribbon.Cap, Kauffman): pass
class Sum(ribbon.Sum, Kauffman): pass

Kauffman.cup_factory = Cup
Kauffman.cap_factory = Cap
Kauffman.sum_factory = Sum

class Variable(ribbon.Box, Kauffman): pass

Kauffman.braid = lambda x, y: (Variable('A', 0, 0) @ x @ y)\
    + (Cup(x, y) >> Variable('A', 0, 0).dagger() >> Cap(x, y))

K = ribbon.Functor(ob=lambda _: 1, ar={}, cod=Category(ribbon.PRO, Kauffman))
drawing.Equation(link, K(link), symbol="$\\mapsto$").draw()
\end{minted}
\[\input{\tikzitpathfigures/kauffman.tikz}\]
\end{figure}

\begin{figure}[H]
\caption{Using Python function syntax to define a string diagram with copy and discard.}
\label{fig:copy-discard}
\begin{minipage}{0.68\textwidth}
\begin{minted}{python}
from discopy.markov import *

x, y = Ty('x'), Ty('y')
f = Box('f', x @ x, y)

@Diagram.from_callable(x @ x, y)
def diagram(a, b):  # Take two wires as inputs
    _ = f(b, a)     # Swap, apply f and discard the result.
    return f(a, b)  # Apply f again and return the result.

assert diagram == Copy(x) @ Copy(x)\
    >> x @ (Swap(x, x) >> f >> Discard(y)) @ x >> f
\end{minted}
\end{minipage}
\begin{minipage}{0.32\textwidth}
\[\begin{tikzpicture}[baseline={(0.base)}, scale=1.4]
	\begin{pgfonlayer}{nodelayer}
		\node [style=none] (1) at (0.625, 4.4) {};
		\node [style=none] (2) at (0.625, 4.15) {};
		\node [style=none] (3) at (2.375, 4.4) {};
		\node [style=none] (4) at (2.375, 3.525) {};
		\node [style=none] (5) at (0.625, 3.925) {};
		\node [style=none] (6) at (0.25, 3.675) {};
		\node [style=none] (7) at (1, 3.675) {};
		\node [style=none] (8) at (0.25, 0.85) {};
		\node [style=none] (9) at (1, 2.9) {};
		\node [style=none] (10) at (2.375, 3.275) {};
		\node [style=none] (11) at (2, 3.025) {};
		\node [style=none] (12) at (2.75, 3.025) {};
		\node [style=none] (13) at (2, 2.9) {};
		\node [style=none] (14) at (2.75, 0.85) {};
		\node [style=none] (15) at (1.5, 2.65) {};
		\node [style=none] (16) at (1, 2.4) {};
		\node [style=none] (17) at (2, 2.4) {};
		\node [style=none] (18) at (1, 2.15) {};
		\node [style=none] (19) at (2, 2.15) {};
		\node [style=none] (20) at (1.5, 1.65) {};
		\node [style=none] (21) at (1.5, 1.4) {};
		\node [style=none] (22) at (1.5, 1.275) {};
		\node [style=none] (23) at (1.5, 0.35) {};
		\node [style=none] (24) at (1.5, 0) {};
		\node [circle, fill=black, scale=0.408] (25) at (0.625, 3.925) {};
		\node [circle, fill=black, scale=0.408] (25) at (2.375, 3.275) {};
		\node [circle, fill=black, scale=0.408] (25) at (1.5, 1.275) {};
		\node [style=none] (25) at (0.75, 1.65) {};
		\node [style=none] (26) at (2.25, 1.65) {};
		\node [style=none] (27) at (2.25, 2.15) {};
		\node [style=none] (28) at (0.75, 2.15) {};
		\node [style=none, fill=white] (29) at (1.5, 1.9) {\py{f}};
		\node [style=none] (30) at (0, 0.35) {};
		\node [style=none] (31) at (3, 0.35) {};
		\node [style=none] (32) at (3, 0.85) {};
		\node [style=none] (33) at (0, 0.85) {};
		\node [style=none, fill=white] (34) at (1.5, 0.6) {\py{f}};
	\end{pgfonlayer}
	\begin{pgfonlayer}{edgelayer}
		\draw [in=90, out=-90] (1.center) to (2.center);
		\draw [in=90, out=-90] (3.center) to (4.center);
		\draw [in=90, out=180, looseness=0.94] (5.center) to (6.center);
		\draw [in=90, out=0, looseness=0.94] (5.center) to (7.center);
		\draw [in=90, out=-90] (2.center) to (5.center);
		\draw [in=90, out=-90] (6.center) to (8.center);
		\draw [in=90, out=-90] (7.center) to (9.center);
		\draw [in=90, out=180, looseness=0.94] (10.center) to (11.center);
		\draw [in=90, out=0, looseness=0.94] (10.center) to (12.center);
		\draw [in=90, out=-90] (4.center) to (10.center);
		\draw [in=90, out=-90] (11.center) to (13.center);
		\draw [in=90, out=-90] (12.center) to (14.center);
		\draw [in=90, out=180, looseness=0.94] (15.center) to (16.center);
		\draw [in=90, out=0, looseness=0.94] (15.center) to (17.center);
		\draw [in=180, out=-90, looseness=0.94] (9.center) to (15.center);
		\draw [in=0, out=-90, looseness=0.94] (13.center) to (15.center);
		\draw [in=90, out=-90] (16.center) to (18.center);
		\draw [in=90, out=-90] (17.center) to (19.center);
		\draw [in=90, out=-90] (20.center) to (21.center);
		\draw [in=90, out=-90] (21.center) to (22.center);
		\draw [in=90, out=-90] (23.center) to (24.center);
		\draw [-, fill=white] (25.center)
			 to (26.center)
			 to (27.center)
			 to (28.center)
			 to cycle;
		\draw [-, fill=white] (30.center)
			 to (31.center)
			 to (32.center)
			 to (33.center)
			 to cycle;
	\end{pgfonlayer}
\end{tikzpicture}\]
\end{minipage}
\end{figure}

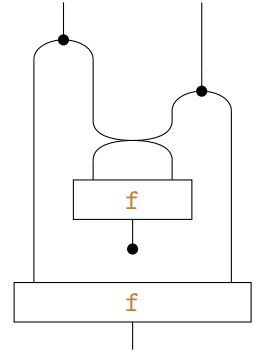
\begin{figure}[H]
\caption{The hypergraph representation of a diagram with spiders, a.k.a. Frobenius algebras.}
\label{fig:hygergraph}
\begin{minted}{python}
from discopy.frobenius import *

x, y = Ty('x'), Ty('y')
f, g = Box('f', x, y), Box('g', y @ y, x)

diagram_lhs = Swap(y, x) >> x @ Cap(x, x) @ y >> Spider(2, 2, x) @ f @ y >> x @ x @ g\
    >> x @ Cup(x, x) @ Spider(0, 0, x)

diagram_rhs = Cap(y, y) @ y @ x >> y @ g @ x >> y @ Spider(2, 2, x) @ Cap(x, x)\
    >> y @ f @ x @ Cup(x, x) >> Cup(y, y) @ x

a, b, c, d = "abcd"
hypergraph = Hypergraph(
    dom=y @ x, cod=x, boxes=(f, g),
    wires=((c, a),           # input wires of the hypergraph
           (((a, ), (b, )),    # input and output wires of f
            ((b, c), (a, ))),  # input and output wires of g
           (a, )),           # output wire of the hypergraph
    spider_types={a: x, b: y, c: y, d: x})  # note the extra x

assert diagram_lhs.to_hypergraph() == hypergraph == diagram_rhs.to_hypergraph()
\end{minted}
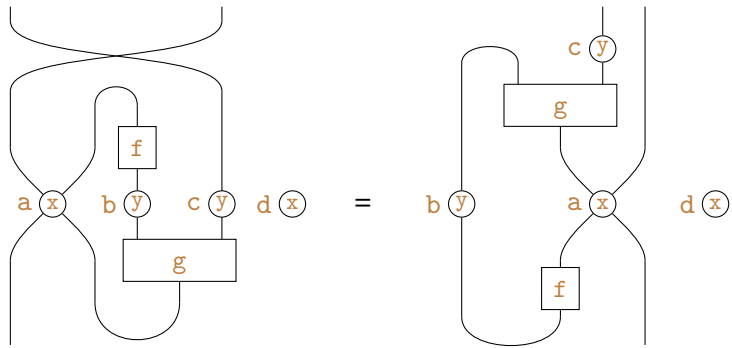
\[\begin{tikzpicture}
	\begin{pgfonlayer}{nodelayer}
		\node [style=none] (20) at (-0.25, 5.75) {};
		\node [style=none] (21) at (3.5, 5.75) {};
		\node [style=none] (22) at (3.5, 4.5) {};
		\node [style=none] (23) at (-0.25, 4.5) {};
		\node [style=none] (24) at (-0.25, 6) {};
		\node [style=none] (25) at (3.5, 6) {};
		\node [style=none] (28) at (-0.25, 3.5) {};
		\node [style=none] (29) at (1.25, 3.5) {};
		\node [style=none] (30) at (-0.25, 1.5) {};
		\node [style=vertex set, minimum width=.35cm] (31) at (0.5, 2.5) {\py{x}};
		\node [style=none] (32) at (1.25, 1.5) {};
		\node [style=none] (33) at (0, 2.5) {\py{a}};
		\node [style=none] (34) at (2.75, 0.75) {};
		\node [style=gate, minimum width=1.5cm] (35) at (2.75, 1.5) {\py{g}};
		\node [style=none] (38) at (2, 1.5) {};
		\node [style=none] (39) at (3.5, 1.5) {};
		\node [style=none] (40) at (2, 4.25) {};
		\node [style=gate, minimum width=.5cm] (41) at (2, 3.5) {\py{f}};
		\node [style=vertex set, minimum width=.35cm] (43) at (2, 2.5) {\py{y}};
		\node [style=none] (44) at (1.5, 2.5) {\py{b}};
		\node [style=none] (45) at (1.25, 4.25) {};
		\node [style=none] (46) at (1.25, 0.75) {};
		\node [style=none] (47) at (-0.25, 0) {};
		\node [style=vertex set, minimum width=.35cm] (48) at (3.5, 2.5) {\py{y}};
		\node [style=none] (49) at (3, 2.5) {\py{c}};
		\node [style=vertex set, minimum width=.35cm] (50) at (4.75, 2.5) {\py{x}};
		\node [style=none] (51) at (4.25, 2.5) {\py{d}};
		\node [style=none] (52) at (6, 2.5) {=};
		\node [style=none] (57) at (10.25, 6) {};
		\node [style=none] (58) at (11, 6) {};
		\node [style=none] (59) at (11, 3.5) {};
		\node [style=none] (60) at (9.5, 3.5) {};
		\node [style=none] (61) at (11, 1.5) {};
		\node [style=vertex set, minimum width=.35cm] (62) at (10.25, 2.5) {\py{x}};
		\node [style=none] (63) at (9.5, 1.5) {};
		\node [style=none] (64) at (9.75, 2.5) {\py{a}};
		\node [style=none] (65) at (9.5, 3.5) {};
		\node [style=gate, minimum width=1.5cm] (66) at (9.5, 4.25) {\py{g}};
		\node [style=none] (67) at (8.75, 4.25) {};
		\node [style=none] (68) at (10.25, 4.25) {};
		\node [style=none] (69) at (9.5, 0.5) {};
		\node [style=gate, minimum width=.5cm] (70) at (9.5, 1) {\py{f}};
		\node [style=vertex set, minimum width=.35cm] (71) at (7.75, 2.5) {\py{y}};
		\node [style=none] (72) at (7.25, 2.5) {\py{b}};
		\node [style=none] (75) at (11, 0) {};
		\node [style=vertex set, minimum width=.35cm] (76) at (10.25, 5.25) {\py{y}};
		\node [style=none] (77) at (9.75, 5.25) {\py{c}};
		\node [style=none] (79) at (8.75, 5) {};
		\node [style=none] (80) at (7.75, 5) {};
		\node [style=none] (81) at (7.75, 0.5) {};
		\node [style=vertex set, minimum width=.35cm] (82) at (12.25, 2.5) {\py{x}};
		\node [style=none] (83) at (11.75, 2.5) {\py{d}};
	\end{pgfonlayer}
	\begin{pgfonlayer}{edgelayer}
		\draw [in=-90, out=90, looseness=0.50] (23.center) to (21.center);
		\draw [in=-90, out=90, looseness=0.50] (22.center) to (20.center);
		\draw (24.center) to (20.center);
		\draw (25.center) to (21.center);
		\draw [in=90, out=-135, looseness=0.75] (31) to (30.center);
		\draw [in=-45, out=90, looseness=0.75] (32.center) to (31);
		\draw [in=-90, out=135, looseness=0.75] (31) to (28.center);
		\draw [in=45, out=-90, looseness=0.75] (29.center) to (31);
		\draw (23.center) to (28.center);
		\draw (34.center) to (35);
		\draw (41) to (40.center);
		\draw (38.center) to (43);
		\draw (22.center) to (39.center);
		\draw (41) to (43);
		\draw [bend left=270, looseness=1.50] (40.center) to (45.center);
		\draw (29.center) to (45.center);
		\draw [bend right=90, looseness=1.50] (46.center) to (34.center);
		\draw (46.center) to (32.center);
		\draw (47.center) to (30.center);
		\draw [in=90, out=-45, looseness=0.75] (62) to (61.center);
		\draw [in=-135, out=90, looseness=0.75] (63.center) to (62);
		\draw [in=-90, out=45, looseness=0.75] (62) to (59.center);
		\draw [in=135, out=-90, looseness=0.75] (60.center) to (62);
		\draw (65.center) to (66);
		\draw (70) to (69.center);
		\draw (75.center) to (61.center);
		\draw (57.center) to (76);
		\draw (76) to (68.center);
		\draw (58.center) to (59.center);
		\draw (79.center) to (67.center);
		\draw [bend left=90] (80.center) to (79.center);
		\draw (71) to (80.center);
		\draw (70) to (63.center);
		\draw [bend right=90] (81.center) to (69.center);
		\draw (71) to (81.center);
	\end{pgfonlayer}
\end{tikzpicture}\]
\end{figure}

\begin{figure}[H]
\caption{The first-order logic formula
$\exists x \cdot G(x) \land \forall y \cdot (\neg G(y) \lor x = y) \land \exists z \cdot M(z) \land P(z, x)$ as a diagram with boxes as predicates, spiders as variables and bubbles as negation, as pioneered by Peirce~\cite{Peirce06}.}
\label{fig:peirce}
\begin{minipage}{0.67\textwidth}
\begin{minted}{python}
from discopy.frobenius import *
from discopy.tensor import Dim, Tensor

Tensor[bool].bubble = lambda self, **_: self.map(lambda x: not x)

@factory
class Formula(Diagram):
    ty_factory = PRO

    def eval(self, size, model):
        return Functor(
            ob=lambda _: Dim(size), ar=lambda f: model[f],
            cod=Category(Dim, Tensor[bool]))(self)

class Cut(Bubble, Formula): pass
class Ligature(Spider, Formula): pass
class Predicate(Box, Formula): pass

P = Predicate("P", 0, 2)  # A binary predicate, i.e. a relation.
G, M = [Predicate(unary, 0, 1) for unary in ("G", "M")]
p, g, m = [[0, 1], [0, 0]], [0, 1], [1, 0]
size, model = 2, {G: g, M: m, P: p}

formula = G >> Ligature(1, 2, PRO(1))\
    >> Cut(Cut(Formula.id(1)) >> G.dagger())\
    @ (M @ 1 >> P.dagger())

assert bool(formula.eval(size, model)) == any(
    g[x] and all(not g[y] or x == y for y in range(size))
    and m[z] and p[z][x] for x in range(size) for z in range(size))
\end{minted}
\end{minipage}
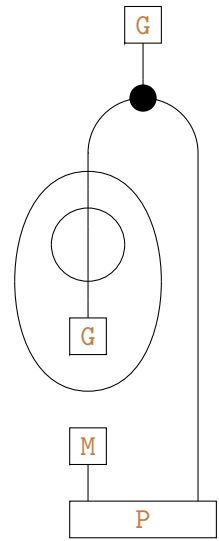
\begin{minipage}{0.33\textwidth}
\[\begin{tikzpicture}[baseline={(0.base)}, scale=1.3]
	\begin{pgfonlayer}{nodelayer}
		\node [style=none] (1) at (2.75, 7.5) {};
		\node [style=none] (2) at (2.75, 7) {};
		\node [style=none] (4) at (2.75, 6.75) {};
		\node [style=none] (5) at (2, 6) {};
		\node [style=none] (6) at (3.5, 6) {};
		\node [style=none] (7) at (2, 6) {};
		\node [style=none] (9) at (3.5, 1.25) {};
		\node [style=none] (11) at (2, 5.75) {};
		\node [style=none] (12) at (1, 4.25) {};
		\node [style=none] (13) at (3, 4.25) {};
		\node [style=none] (14) at (2, 5.75) {};
		\node [style=none] (15) at (1, 4.25) {};
		\node [style=none] (17) at (2, 5.5) {};
		\node [style=none, fill=white, right] (18) at (2.1, 5.65) {};
		\node [style=none] (19) at (3, 4.25) {};
		\node [style=none] (21) at (2, 5.25) {};
		\node [style=none] (22) at (1.5, 4.75) {};
		\node [style=none] (23) at (2.5, 4.75) {};
		\node [style=none] (24) at (2, 5) {};
		\node [style=none] (25) at (1.5, 4.75) {};
		\node [style=none] (27) at (2, 4.5) {};
		\node [style=none, fill=white, right] (28) at (2.1, 4.9) {};
		\node [style=none] (29) at (2.5, 4.75) {};
		\node [style=none] (31) at (2, 4.25) {};
		\node [style=none] (32) at (2, 4) {};
		\node [style=none] (33) at (2, 3.75) {};
		\node [style=none, fill=white, right] (34) at (2.1, 3.9) {};
		\node [style=none] (35) at (2, 2.75) {};
		\node [style=none] (36) at (2, 1.75) {};
		\node [style=none] (37) at (2, 1.25) {};
		\node [circle, fill=black] (39) at (2.75, 6.75) {};
		\node [style=none] (39) at (2.5, 7.5) {};
		\node [style=none] (40) at (3, 7.5) {};
		\node [style=none] (41) at (3, 8) {};
		\node [style=none] (42) at (2.5, 8) {};
		\node [style=none, fill=white] (43) at (2.75, 7.75) {\py{G}};
		\node [style=none] (44) at (1.75, 3.25) {};
		\node [style=none] (45) at (2.25, 3.25) {};
		\node [style=none] (46) at (2.25, 3.75) {};
		\node [style=none] (47) at (1.75, 3.75) {};
		\node [style=none, fill=white] (48) at (2, 3.5) {\py{G}};
		\node [style=none] (49) at (1.75, 1.75) {};
		\node [style=none] (50) at (2.25, 1.75) {};
		\node [style=none] (51) at (2.25, 2.25) {};
		\node [style=none] (52) at (1.75, 2.25) {};
		\node [style=none, fill=white] (53) at (2, 2) {\py{M}};
		\node [style=none] (54) at (1.75, 0.75) {};
		\node [style=none] (55) at (3.75, 0.75) {};
		\node [style=none] (56) at (3.75, 1.25) {};
		\node [style=none] (57) at (1.75, 1.25) {};
		\node [style=none, fill=white] (58) at (2.75, 1) {\py{P}};
	\end{pgfonlayer}
	\begin{pgfonlayer}{edgelayer}
		\draw [in=90, out=-90] (1.center) to (2.center);
		\draw [in=90, out=-180] (4.center) to (5.center);
		\draw [in=90, out=0] (4.center) to (6.center);
		\draw [in=90, out=-90] (2.center) to (4.center);
		\draw [in=90, out=-90] (5.center) to (7.center);
		\draw [in=90, out=-90] (6.center) to (9.center);
		\draw [in=90, out=180] (11.center) to (12.center);
		\draw [in=90, out=0] (11.center) to (13.center);
		\draw [in=90, out=-90] (7.center) to (14.center);
		\draw [in=90, out=-90] (14.center) to (17.center);
		\draw [in=90, out=180] (21.center) to (22.center);
		\draw [in=90, out=0] (21.center) to (23.center);
		\draw [in=90, out=-90] (17.center) to (24.center);
		\draw [in=90, out=-90] (24.center) to (27.center);
		\draw [in=-180, out=-90] (25.center) to (31.center);
		\draw [in=90, out=-90] (27.center) to (32.center);
		\draw [in=0, out=-90] (29.center) to (31.center);
		\draw [in=90, out=-90] (32.center) to (33.center);
		\draw [in=-180, out=-90] (15.center) to (35.center);
		\draw [in=0, out=-90] (19.center) to (35.center);
		\draw [in=90, out=-90] (36.center) to (37.center);
		\draw [-, fill=white] (39.center)
			 to (40.center)
			 to (41.center)
			 to (42.center)
			 to cycle;
		\draw [-, fill=white] (44.center)
			 to (45.center)
			 to (46.center)
			 to (47.center)
			 to cycle;
		\draw [-, fill=white] (49.center)
			 to (50.center)
			 to (51.center)
			 to (52.center)
			 to cycle;
		\draw [-, fill=white] (54.center)
			 to (55.center)
			 to (56.center)
			 to (57.center)
			 to cycle;
	\end{pgfonlayer}
\end{tikzpicture}\]
\end{minipage}
\end{figure}

\appendix

\section{Free categories and functors in Python}\label{cat}

The most basic module of DisCoPy is \pyref{cat}{cat},\footnote
{To make this report easy to use, the first mention of each module and class is a clickable link to the documentation.} an implementation of the \emph{free category}\footnote
{More precisely, what we mean is the free category generated by the signature with \py{Ob} as objects and \py{Box} as arrows.}
with the class \pyref{cat.Ob}{Ob}\py{(name: str)} as objects and the class \pyref{cat.Arrow}{Arrow}{}\py{(inside: list[Arrow], dom: Ob, cod: Ob)} as arrows.\footnote{For technical reasons, the implementation of DisCoPy uses the immutable \py{tuple[X, ...]} rather than \py{list[X]}.}
That is, an arrow \py{f} is encoded as a list of arrows \py{f.inside} with a domain and a codomain.
The method \py{Arrow.id(dom: Ob) -> Arrow} returns an arrow with an empty list \py{inside} while the method \py{Arrow.then(self, *others: Arrow) -> Arrow} does concatenation or raises \py{AxiomError} if the arrows do not compose.
They are shortened to \py{Id} and the binary operator \py{>>} respectively.

\pyref{cat.Box}{Box}{}\py{(name: str, dom: Ob, cod: Ob)} implements generating arrows with \py{f.inside = [f]}
so that we have \py{f >> Id(f.cod) == f == Id(f.dom) >> f} on the nose.
Boxes have an optional attribute \py{data: Any} which can be used to parameterise arrows with SymPy~\cite{MeurerEtAl17} symbols.
This comes with a property \py{Arrow.free_symbols} and two methods \py{Arrow.subs} for substitution and \py{Arrow.lambdify} for fast numerical computation.
An optional attribute \py{is_dagger: bool} implements \emph{free dagger categories} with the method \py{Arrow.dagger(self) -> Arrow}, shortened to the list-reversal operator \py{[::-1]}.

\pyref{cat.Sum}{Sum}{}\py{(terms: list[Arrow], dom: Ob, cod: Ob)} is a subclass of \py{Box} which implements \emph{enrichment over commutative monoids}, i.e. formal sums of parallel arrows with the method \py{Arrow.zero (dom: Ob, cod: Ob) -> Sum} as unit.
Composition of sums is implemented so that the equations \py{f >> (g + g_) == f >> g + f >> g_} and \py{(f + f_) >> g == f >> g + f_ >> g} hold on the nose.
\pyref{cat.Bubble}{Bubble}{}\py{(arg: Arrow, dom: Ob, cod: Ob)} is a subclass of \py{Box} which allows to encode \emph{unary operators on homsets}, i.e. arbitrary functions from arrows to arrows.

\pyref{cat.Category}{Category}{}\py{(ob: type, ar: type)} is just a pair of types for objects and arrows with methods \py{dom, cod, id, then} of the appropriate type.
For instance, \py{Pyth = Category(type, Function)} has Python types as objects and \pyref{python.Function}{Function}{}\py{(inside: Callable, dom: type, cod: type)} as arrows;
\py{Mat[R] = Category(int, Matrix[R])} has natural numbers as objects and \pyref{matrix.Matrix}{Matrix}{}\py{[R](inside: array, dom: int, cod: int)} as arrows with \py{array} from any of NumPy~\cite{VanDerWaltEtAl11}, TensorFlow~\cite{AbadiEtAl16}, PyTorch~\cite{PaszkeEtAl19} or JAX~\cite{BradburyEtAl18} and entries in any rig \py{R}, i.e. any data type with sum, product, zero and one.

\pyref{cat.Functor}{Functor}{}\py{(ob: Map, ar: Map, cod: Category)} is given by an optional codomain and a pair of \py{Map = dict | Callable} from \py{Ob} to \py{cod.ob} and from \py{Box} to \py{cod.ar}.
By default, the codomain is the free \py{Category(Ob, Arrow)}.
The domain is defined implicitly by the domain of \py{ob} and \py{ar}.
Functors also have their own methods \py{id} and \py{then} so we can define \py{CAT = Category(Category, Functor)}.

\section{Planar diagrams for monoidal categories}\label{monoidal}

The \pyref{monoidal}{monoidal} module is where planar diagrams are implemented.
\pyref{monoidal.Ty}{Ty}{}\py{(inside: list[Ob])} is a subclass of \py{Ob} with a method \py{Ty.tensor(self, *others: Ty) -> Ty} for concatenation shortened to the binary operator \py{@} with the empty type \py{Ty()} as unit, i.e. \py{Ty} is the free monoid over \py{Ob}.

\pyref{monoidal.Layer}{Layer}{}\py{(left: Ty, box: Box, right: Ty, *more: Ty | Box)} is a subclass of \py{Box} with two methods for \emph{whiskering}, i.e. concatenating a layer \py{f} with a type \py{x} on the left \py{x @ f} and right \py{f @ x}.
It also comes with a method \py{Layer.tensor(self, *others: Layer) -> Layer} which returns a layer with potentially many boxes alternating with types between them (this is a new feature of DisCoPy v1.0).

\pyref{monoidal.Diagram}{Diagram}{}\py{(inside: list[Layer], dom: Ty, cod: Ty)} is a subclass of \py{Arrow} with layers as boxes and a method \py{Diagram.tensor(self, *others: Diagram) -> Diagram} shortened to \py{@} and defined in terms of whiskering and composition so that \py{f @ g = f @ g.dom >> f.cod @ g}.

The method \py{Diagram.draw} plots the diagram with either Matplotlib~\cite{Hunter07} or TikZ~\cite{Tantau13}.
The method \py{Diagram.normal_form} solves the word problem for (connected) planar diagrams~\cite{DelpeuchVicary22} by repeatedly applying \emph{interchanger} rewrites from \py{f.dom @ g >> f @ g.cod} to \py{f @ g.dom >> f.cod @ g}.
Note that by default, two diagrams are equal only if their tuple of layers \py{inside} are, i.e. DisCoPy implements the \emph{free premonoidal category}~\cite{PowerRobinson97}.
This is a feature rather than a bug, indeed when functions have side-effects \py{Pyth} is only premonoidal.

\pyref{monoidal.Box}{Box}{}\py{(name: str, dom: Ty, cod: Ty)} is a subclass of \py{cat.Box} and \py{Diagram} with \py{f.inside = [Layer(Ty(), f, Ty())]} so that we have \py{f @ Id() == f == Id() @ f} on the nose, where \py{Id()} is the empty diagram, i.e. the identity on the empty type.
\pyref{monoidal.Functor}{monoidal.Functor} is a subclass of \pyref{cat.Functor}{cat.Functor} with \py{Category(Ty, Diagram)} as domain and an arbitrary monoidal category as codomain.
This includes \py{Pyth} with \py{tuple} as tensor, \py{Mat[R]} with direct sum as well as \py{Category(Dim, Tensor)} where \py{Dim} is the free monoid over positive integers and \py{Tensor} is a subclass of \py{Matrix} with the Kronecker product as tensor.
Monoidal functors into \py{Tensor} correspond to \emph{tensor network contraction}, which DisCoPy computes via the specialised TensorNetwork library~\cite{RobertsEtAl19}.

DisCoPy then goes on to implement the hierarchy of graphical languages for monoidal categories as described in Selinger's survey~\cite{Selinger10}, see Figure~\ref{hierarchy}.
Structural morphisms for types of length one are implemented as subclasses of \py{Box}, e.g. \pyref{braided.Braid}{Braid}, while the structural morphisms for arbitrary types are implemented as methods, e.g. \py{Diagram.braid}, which ensure that coherence laws hold on the nose.
The functor class in each module respects the structure, so that e.g. \py{braided.Functor} sends each \py{Braid} box in its domain to the \py{braid} method of its codomain category, see Figure~\ref{fig:kauffman}.

A notable addition of DisCoPy v1.0 is the implementation of \emph{free traced categories}~\cite{JoyalEtAl96} where the trace can be interpreted as partial matrix trace in \py{Tensor}, as reflexive transitive closure in \py{Mat[bool]} or as parameterised fixed points in \py{Pyth} as demonstrated in Figure~\ref{traced-example}.
This comes with an \pyref{interaction}{interaction} module for the \emph{Int-construction}, also called the \emph{geometry of interaction}~\cite{Abramsky96}, which turns any balanced (symmetric) traced category into a free ribbon (compact) category.

Note that the \py{traced} module implements planar traced categories, of which \py{pivotal} is an example.
We avoid extra modules for each kind of traced category, so the \py{balanced} and \py{symmetric} modules come with traces by default.
This is justified by the fact that the free balanced category can be faithfully embedded in the free balanced traced category.

\section{Hypergraphs for symmetric categories}\label{hypergraph}

The encoding of string diagrams as lists of layers makes it possible to draw them and evaluate them simply with a \py{for} loop.
However this comes at the cost of representing swaps explicitly as generating morphisms subject to naturality conditions.
Alternatively, string diagrams for symmetric categories can be encoded as
\emph{discrete cospans of hypergraphs}~\cite{BonchiEtAl22} where the equations for symmetric, traced, compact and hypergraph categories all come for free.
This is implemented in DisCoPy's \pyref{hypergraph}{hypergraph} module.

The class \pyref{hypergraph.Hypergraph}{Hypergraph}{}\py{[C](dom, cod, boxes, wires, spider_types)} is defined by:
\begin{itemize}
    \item a type parameter \py{C: Category} with \py{C.ob = Ty} and \py{C.ar = Diagram} by default,
    \item a pair of types \py{dom: C.ob} and \py{cod: C.ob} together with a list of \py{boxes: list[C.ar]},
    \item a mapping \py{spider_types: dict[Spider, C.ob]} with keys of any type \py{Spider = Any},
    \item a triple \py{wires: tuple[W, list[tuple[W, W]], W]} where \py{W = list[Spider]}.
\end{itemize}
The first and last lists of spider \py{wires} correspond to the input and output wires of the overall diagram, while the middle list corresponds to the input and output wires of each box.
That is, we require that:
\begin{itemize}
    \item \py{len(wires[0]) == len(dom)} and \py{len(wires[2]) == len(cod)},
    \item for \py{i, box in enumerate(boxes)} and \py{box_dom_wires, box_cod_wires = wires[1][i]},
    \py{len(box_dom_wires) == len(box.dom)} and \py{len(box_cod_wires) == len(box.cod)}.
\end{itemize}

The method \py{Hypergraph.tensor} concatenates the attributes of two hypergraphs then reorders the wires.
The composition \py{Hypergraph.then} computes the push-out of cospans via reflexive transitive closure.
The three methods \py{Hypergraph.id}, \py{swap} and \py{spiders} generate all the hypergraphs with no boxes.
\py{Hypergraph.from_box} wraps up a box as a hypergraph
while \py{to_diagram} represents the hypergraph as a planar \pyref{frobenius.Diagram}{frobenius.Diagram}, with explicit \pyref{frobenius.Swap}{Swap} and \pyref{frovenius.Spider}{Spider} boxes.
The inverse translation \py{Diagram.to_hypergraph} is a functor with \py{cod=Category(Ty, Hypergraph)}.
Two hypergraphs are equal when their attributes are equal up to a permutation of the boxes and spiders, this is computed by reduction to the graph isomorphism algorithm of NetworkX~\cite{HagbergEtAl08}.
Hypergraphs are the arrows of \emph{free hypergraph categories}, i.e. symmetric categories with a supply of spiders, also known as special commutative Frobenius algebras.

The property \py{Hypergraph.is_bijective: bool} checks if the wires define a bijection between ports, i.e. each spider is connected to either zero or two ports.
Bijective hypergraphs are the arrows of \emph{free compact categories}, their translation to diagram only involves \pyref{compact.Swap}{Swap}, \pyref{compact.Cup}{Cup} and \pyref{compact.Cap}{Cap}.
The property \py{Hypergraph.is_monogamous: bool} checks if furthermore the bijection goes from output ports (i.e. either the domain of the hypergraph or the codomain of a box) to input ports (i.e. either the domain of a box or the codomain of the hypergraph) in which case the diagrams only involve \pyref{traced.Swap}{Swap} and \pyref{traced.Trace}{Trace}.

A third property \py{Hypergraph.is_causal: bool} checks if each spider is connected to exactly one output port and to zero or more input ports that all have higher indices in the list.
In this case, the diagrams only involve \pyref{markov.Swap}{Swap}, \pyref{markov.Copy}{Copy} and \pyref{markov.Discard}{Discard} boxes which are defined in the \pyref{markov}{markov} module.
Causal hypergraphs are the arrows of the free symmetric category with a supply of cocommutative comonoids, also called a \emph{copy-discard category}.
The method \py{markov.Diagram.normal_form}, still under development at the time of writing, repeatedly applies the equation \py{f >> Discard() == f} in order to enforce the monoidal unit \py{Ty()} as a terminal object.
The equivalence classes of causal hypergraphs are the arrows of the free \emph{Markov category}~\cite{FritzLiang23}.

We also define a weaker property \py{Hypergraph.is_left_monogamous: bool}, which checks if the wires define a function from input ports to output ports, i.e. each spider is either disconnected (i.e. the trace of an identity wire) or connected to exactly one output port.
Left monogamous hypergraphs are the arrows of the \emph{free copy-discard traced category} where diagrams involve \pyref{markov.Swap}{Swap}, \pyref{markov.Copy}{Copy}, \pyref{markov.Discard}{Discard} and \pyref{markov.Trace}{Trace}.
Finally, monogamous causal hypergraphs are the arrows of \emph{free symmetric categories}, their translation to diagrams only involves \pyref{symmetric.Swap}{Swap}.

A powerful new feature built on top of the \py{Hypergraph} class allows to construct string diagrams using the standard syntax for Python functions, a form of \emph{substructural type system} implemented as a decorator (i.e. a higher-order function) which is illustrated in Figure~\ref{fig:copy-discard}.
In practice, this allows the user to easily define morphisms in any copy-discard category where the swapping, copying and discarding of Python variables is encoded explicitly.
It also allows the definition of morphisms in a symmetric category where copying and discarding of variables leads to a type error (e.g. in a quantum circuit) or in a monoidal category where even reordering variables is forbidden (e.g. in a grammatical derivation).

\bibliographystyle{eptcs}
\hbadness=99999
\bibliography{act2023-discopy}
\vfill
\end{document}